\RequirePackage{fix-cm}
\documentclass{svjour3}                     
\smartqed  
\usepackage{graphicx}
\usepackage{amssymb}
\usepackage{amsmath}
\usepackage{xcolor}
\usepackage{subcaption}
\usepackage{hyperref}
\usepackage{algorithm}
\usepackage[noend]{algpseudocode}

\usepackage{tikz}
\usetikzlibrary{shapes.geometric, arrows}

\tikzstyle{startstop} = [rectangle, rounded corners, minimum width=3cm, minimum height=0cm,text centered, draw=black, fill=red!30]
\tikzstyle{process} = [rectangle, minimum width=3cm, minimum height=0cm, text centered, draw=black, fill=blue!30]
\tikzstyle{decision} = [diamond, minimum width=0cm, minimum height=0cm, text centered, draw=black, fill=green!30]
\tikzstyle{arrow} = [thick,->,>=stealth]

\newcommand{\mc}[1]{\mathcal{#1}}
\newcommand{\mb}[1]{\mathbb{#1}}
\newcommand{\empi}{\hat{\mathbb{P}}_N}

\newtheorem{assu}{Assumption}

\newtheorem{prop}{Proposition}
\newtheorem{rema}{Remark}

\begin{document}

\title{Comparative Analysis of Two-Stage Distributionally Robust Optimization over 1-Wasserstein and 2-Wasserstein Balls
}

\titlerunning{Two-Stage Distributionally Robust Optimization over Wasserstein Balls}        

\author{Geunyeong Byeon
}


\institute{Geunyeong Byeon \at
              School of Computing and Augmented Intelligence, Arizona State University, Tempe, AZ\\
              \email{geunyeong.byeon@asu.edu}           %
}

\date{Received: date / Accepted: date}

\maketitle

\begin{abstract}
This paper investigates advantages of using 2-Wasserstein ambiguity sets over 1-Wasserstein sets in two-stage distributionally robust optimization with right-hand side uncertainty. We examine the worst-case distributions within 1- and 2-Wasserstein balls under both unrestricted and nonnegative orthant supports, highlighting a pathological behavior arising in 1-Wasserstein balls. 
Closed-form solutions for a single-scenario newsvendor problem illustrate that 2-Wasserstein balls enable more informed decisions. Additionally, a penalty-based dual interpretation suggests that 2-Wasserstein balls may outperform 1-Wasserstein balls across a broader range of Wasserstein radii, even with general support sets.


\end{abstract}

\section{Introduction}
\label{sec:intro}
Two-stage distributionally robust optimization (TSDRO) seeks an optimal here-and-now decision $x \in \mathbb{R}^{n_x}$ that hedges against the worst-case distribution within an ambiguity set of plausible distributions for uncertain parameters $\tilde{\xi} \in \mathbb{R}^k$ supported on $\Xi$. 
This paper focuses on TSDRO with right-hand side uncertainty, where the ambiguity set is an \(r\)-Wasserstein ball of radius \(\epsilon\) centered at an empirical distribution $\empi$ constructed from $N$ samples, denoted as \(\mathcal{P}^{r}_{N,\epsilon}\). 
We compare the performance of resulting here-and-now decisions under \(\mathcal{P}^{1}_{N,\epsilon}\) and \(\mathcal{P}^{2}_{N,\epsilon}\).

\subsection{Problem statement}
We focus on a class of TSDRO problems in which the first-stage problem is modeled as a mixed-integer conic-linear programming (coneLP) problem and the second-stage problem is posed as a coneLP problem:
\begin{equation}
\begin{aligned}
    \min & \ c^\top x +  \sup_{\mb P \in \mc P^{r}_{N,\epsilon}} \mb E_{\mb P}[Z(x, \tilde\xi)]\\
    \mbox{s.t.} 
    & \ x \in \mc X := \{x \in \mc K_x: x_i \in \mb Z, \ \forall i \in \mc I, \ Ax = b\},
  \end{aligned}    \label{prob:TSDRO}
\end{equation}
where 
\begin{equation}
  Z(x,\xi) := \inf_{y \in \mc K_y}\left\{q^\top y : Wy \ge Tx + \xi \right\}\label{prob:second-stage}
  \end{equation}
evaluates the second-stage cost by optimizing a wait-and-see decision $y$ given $x$ and an outcome $\xi$ of the uncertain vector $\tilde \xi$. Positive integers $n_x$, $n_y$, and $k$ respectively denote the dimensions of $x$, $y$, and $\tilde\xi$. The sets $\mc K_x \subseteq \mb R^{n_x}$ and $\mc K_y \subseteq \mb R^{n_y}$ represent Cartesian products of some collections of proper cones (e.g., second-order cones, nonnegative orthants, and vectorized positive semidefinite cones). The set $\mc X$ defines the feasible region of $x$, where a subvector $(x_i)_{i \in \mathcal I}$ may be integer-valued. 
The parameters $c, b, q, A, W$ and $T$ are vectors and matrices with suitable dimensions.

Given \(N\) samples \(\zeta^1, \cdots, \zeta^N \in \Xi\) of \(\tilde \xi\), we let \(\empi\) denote the empirical distribution \(\empi := \frac{1}{N}\sum^N_{i = 1} {\delta}_{\zeta^i}\), where \({\delta}_{\zeta^i}\) denotes the Dirac measure whose unit mass is concentrated on \(\{\zeta^i\}\). The ambiguity set $\mc P^{r}_{N,\epsilon}$ is given by
\[
\mc P^{r}_{N,\epsilon} := \{\mb Q \in \mc P^r(\Xi): W^{r}(\mb Q, \empi) \le {\epsilon}\},
\] 
where \({\mc P^r(\Xi)}\) denotes the collection of all probability distributions \(\mb Q\) on \((\Xi, \mc B(\Xi))\) with a finite \(r\)th moment, and \(\mc B(\Xi)\) is the Borel \(\sigma\)-algebra of subsets of \(\Xi\). The \(r\)-Wasserstein distance \(W^{r}\) between distributions \(\mb Q_1,\mb Q_2 \in \mc P^r(\Xi)\)
is defined as
\[
W^{r}(\mb Q_1 ,\mb Q_2):=\inf _{\gamma \in \Gamma (\mb Q_1 ,\mb Q_2 )}\left\{\left(\int _{\Xi\times \Xi}\|\xi_1-  \xi_2\|^r\, \mathrm{d}\gamma (\xi_1,\xi_2)\right)^{1/r}\right\},
\]
where $\Gamma (\mb Q_1 ,\mb Q_2 )$ denotes the collection of all probability distributions on the product space $(\Xi \times \Xi, \mc B(\Xi) \otimes \mc B(\Xi))$ with marginals 
$\mb Q_1$ and $\mb Q_2$, respectively, and $\|\cdot\|$ is an $l_p$-norm in $\mb R^k$ for some $p \in [1,\infty]$. 

TSDRO provides an alternative to sample average approximation (SAA) and robust optimization (RO), determined by the Wasserstein radius \(\epsilon > 0\). Let \(\tilde{x}^{r}_{N}(\epsilon)\) denote the random optimal solution of \eqref{prob:TSDRO}, with randomness arising from sampling for $\empi$. When \(\epsilon = 0\), the Wasserstein ambiguity set reduces to \(\mathcal{P}^{r}_{N,0} = \{\empi\}\). In this case, we denote the corresponding solution by \(\tilde{x}^{\texttt{SAA}}_N\). Additionally, let \(\mathcal{X}_{\texttt{RO}} := \arg \inf_{x \in \mathcal{X}} c^\top x + \sup_{\xi \in \Xi} Z(x, \xi)\) denote the set of robust optimization solutions. 
We define $\tilde O^{r}_N(\epsilon)$ as the true expected cost of the solution obtained by TSDRO over $\mathcal P^{r}_{N,\epsilon}$, that is $\tilde O^{r}_N(\epsilon) := c^\top \tilde{x}^{r}_{N}(\epsilon) + \mathbb{E}_{\tilde\xi \sim \mathbb P^{\texttt{true}}}[Z(\tilde{x}^{r}_{N}(\epsilon), \tilde \xi)]$. When $\epsilon = 0$, this value corresponds to the out-of-sample performance of the SAA solution. If ${\mathcal X}_{\texttt{RO}} \neq \emptyset$, then as $\epsilon \rightarrow \infty$, $\tilde O^{r}_N(\epsilon)$ converges to the out-of-sample performance of RO.

\subsection{Motivating example}
\begin{figure}[!t]  
  \centering
  \begin{subfigure}[b]{0.25\textwidth}
      \centering      
      \includegraphics[width=\textwidth]{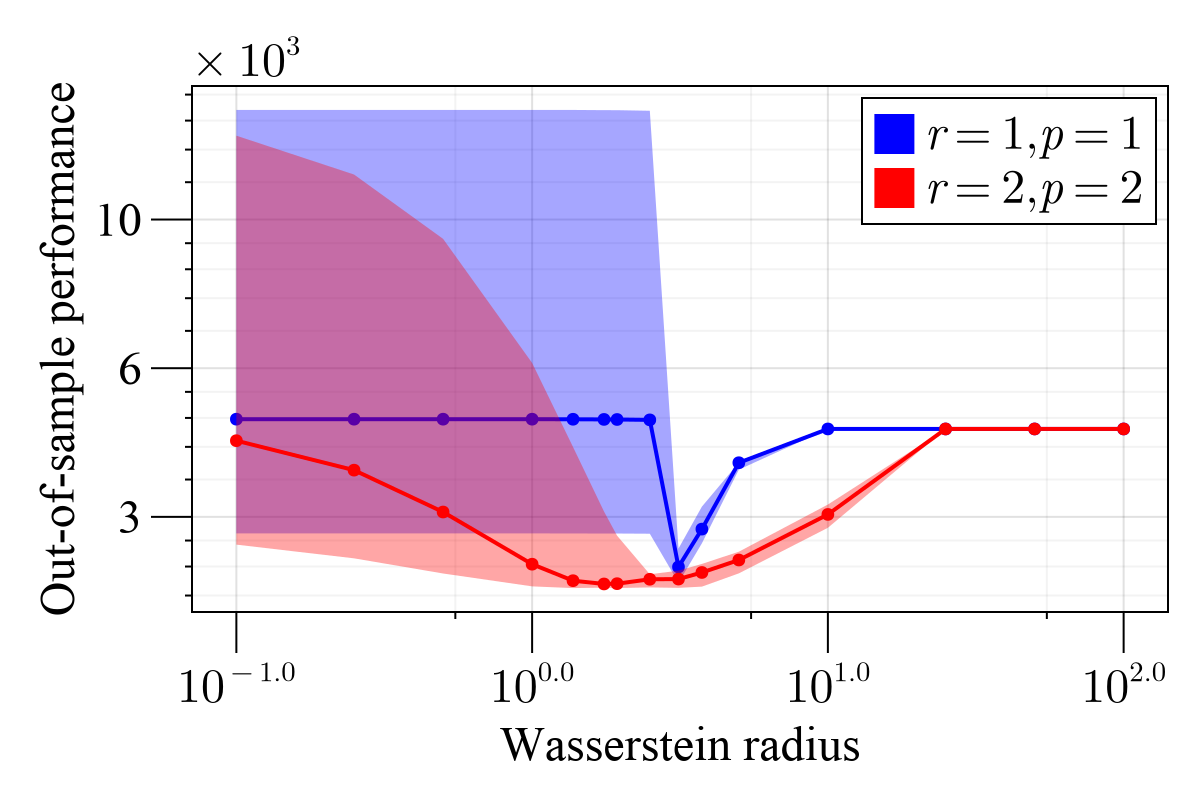}                       
      \caption{\texttt{p1}}\label{fig:oos:p1}
  \end{subfigure}
  \begin{subfigure}[b]{0.25\textwidth}
      \centering      
      \includegraphics[width=\textwidth]{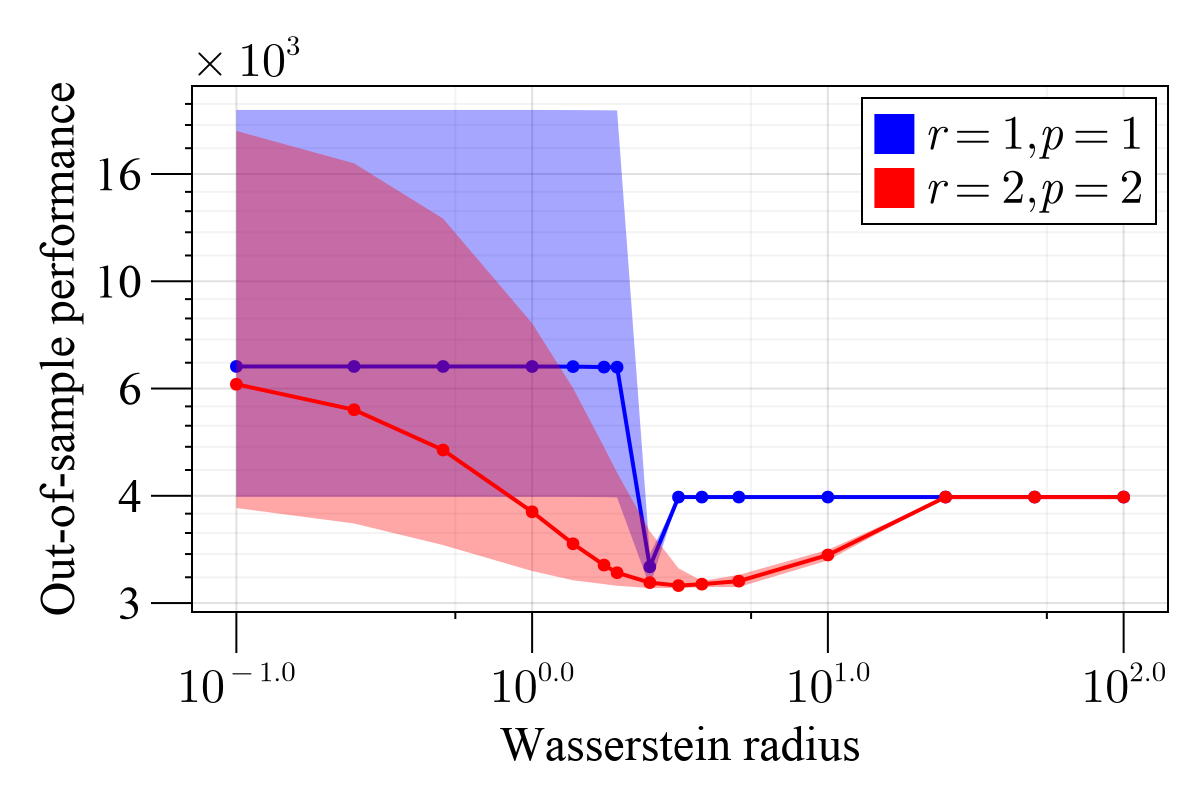}                       
      \caption{\texttt{p21}}\label{fig:oos:p21}
  \end{subfigure} 
  \begin{subfigure}[b]{0.25\textwidth}
      \centering      
      \includegraphics[width=\textwidth]{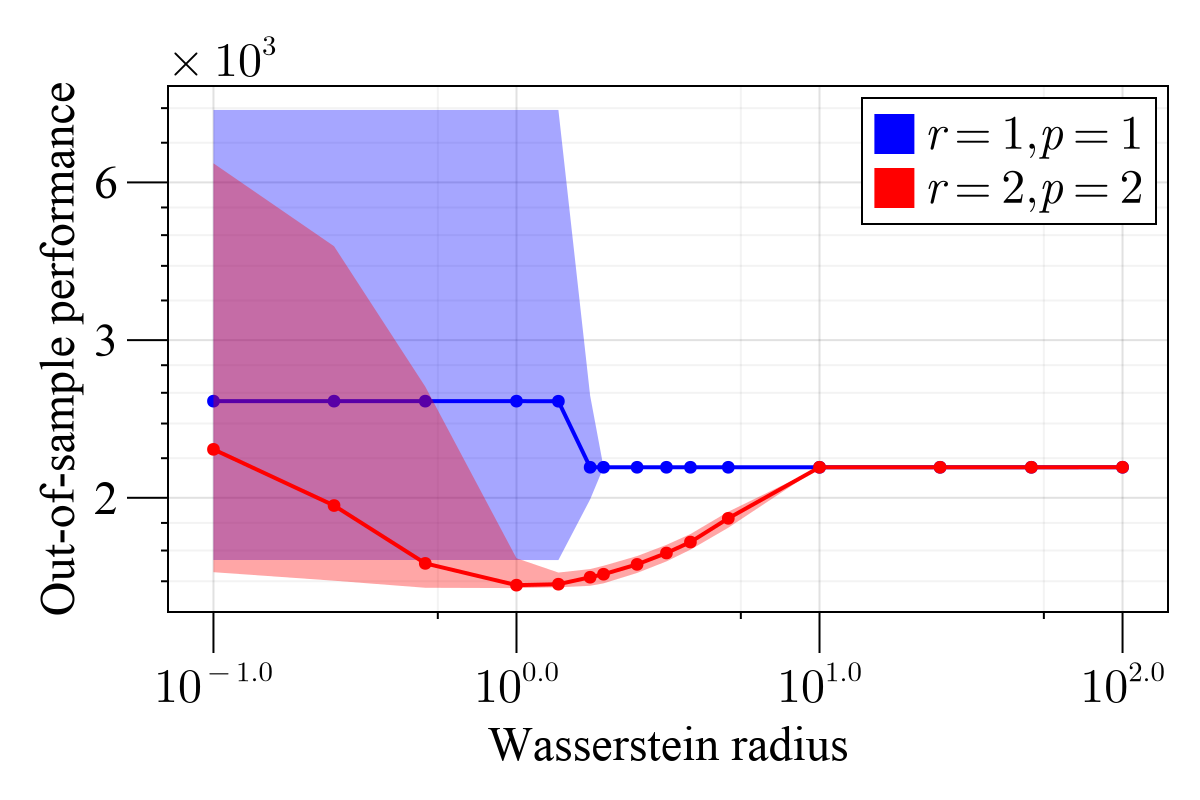}                       
      \caption{\texttt{p43}}\label{fig:oos:p43}
  \end{subfigure} 
    \begin{subfigure}[b]{0.22\textwidth}
      \centering      
      \includegraphics[width=\textwidth]{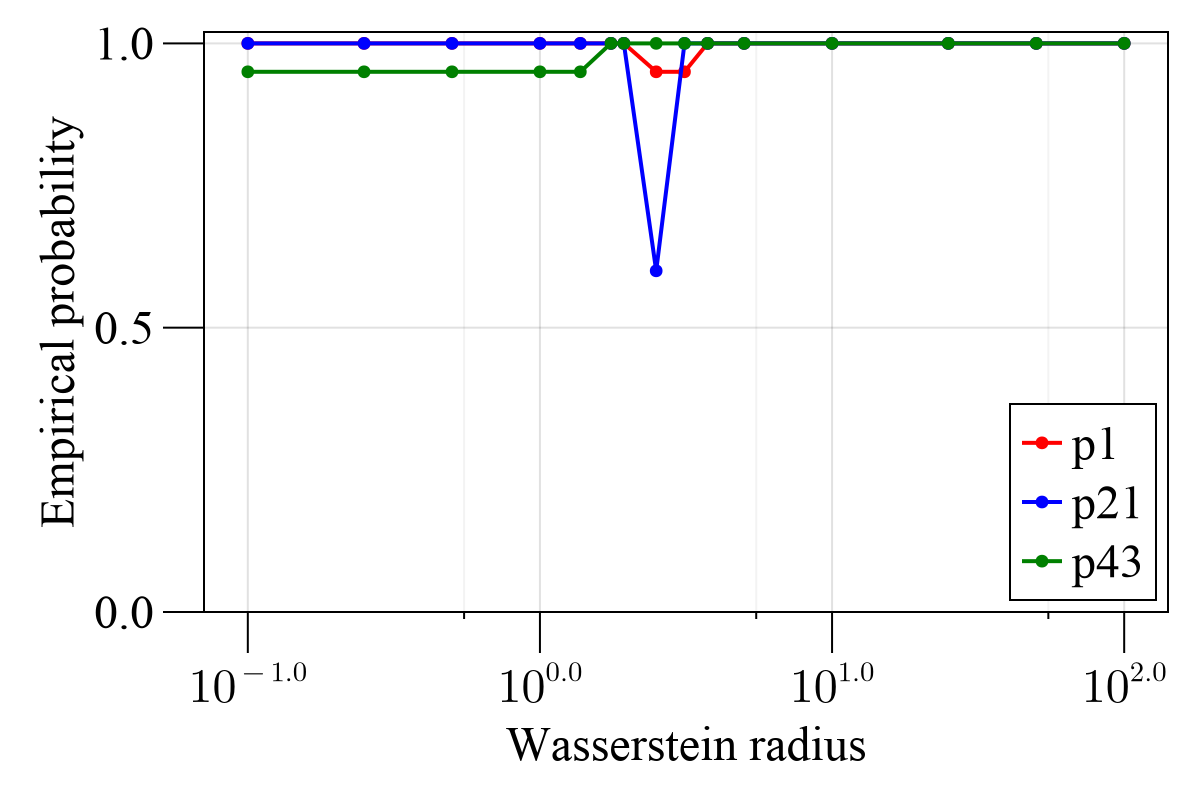}                 
      \caption{
      }
      \label{fig:oos:dominance}
  \end{subfigure} 
  \caption{
  (a)-(c): The median values (line) and the 10th to 90th percentiles (shaded region) of $\tilde O^{r}_5(\epsilon)$,
  where the expectation is estimated using 2,000 testing samples drawn from $\mathbb P^{\texttt{true}}$. (d): $\mathbb P[\tilde O^{2}_5(\epsilon) \le \tilde O^{1}_5(\epsilon) \times 1.001]$. For each choice of $\epsilon$, $\tilde O^{r}_5(\epsilon)$ is simulated over 20 runs} 
  \label{fig:oos}
\end{figure}


 Figure \ref{fig:oos} presents a notable observation from three facility location problem instances—Holmberg instances \texttt{p1}, \texttt{p21}, and \texttt{p43} \cite{holmberg1999exact}—using the experimental setup described in \cite{byeon2022two}, where $p=2$ is used for the $r=2$ case and $p=1$ for the $r=1$ case to ensure comparability on the same scale. The results reveal that $\tilde O^{2}_5(\epsilon)$ has a lower or equal median across all tested 
$\epsilon$ values compared to $\tilde O^{1}_5(\epsilon)$. As shown in the figure, $\tilde x^{1}_5(\epsilon)$ provides alternative solutions to both SAA and RO over a narrow range of $\epsilon$-values. It initially aligns with the SAA solution, then undergoes a sharp decline, and promptly converges to some $x^{\texttt{RO}} \in \mathcal X_{\texttt{RO}}$. In contrast, $\tilde x^2_5(\epsilon)$ exhibits a smoother transition from $\tilde x^{\texttt{SAA}}_5$ as soon as $\epsilon$ is greater than zero and then converges to a RO solution at larger $\epsilon$-values than $\tilde x^{1}_5(\epsilon)$. Furthermore, $\tilde O^{2}_5(\epsilon)$ empirically demonstrates zeroth-order dominance over $\tilde O^{1}_5(\epsilon)$ across a wide range of 
$\epsilon$-values. Specifically, this dominance is observed for values of $\epsilon$ where $\mathbb P[\tilde O^{2}_5(\epsilon) \le \tilde O^{1}_5(\epsilon) \times 1.001]=1$ (see Figure \ref{fig:oos:dominance}). This paper aims to investigate the underlying reasons for this behavior.


\subsection{Related Literature}
TSDRO with right-hand side uncertainty over 1- and 2-Wasserstein ambiguity sets has been extensively studied in the literature \cite{hanasusanto2018conic,esfahani2018data,zhao2018data,gamboa2021decomposition,duque2022distributionally,byeon2022two,mittal2022finding,gao2023distributionally}. However, there is limited research comparing their out-of-sample performances. The most relevant work is \cite{gao2022finite}, which provides performance guarantees for $r = 1$ and $r = 2$. These results, especially for $r=2$, however, rely on the differentiability of $Z(x, \cdot)$, a condition that typically does not hold in our setting.

\subsection{Contributions}
To the best of the author's knowledge, this is the first paper to demonstrate the potential advantages of 2-Wasserstein balls over 1-Wasserstein balls in the TSDRO context. For \( \Xi = \mathbb{R}^k \) or \( \mathbb{R}^k_+ \), we derive the worst-case distributions in \( \mathcal{P}^{1}_{N,\epsilon} \) and \( \mathcal{P}^{2}_{N,\epsilon} \), highlighting an undesirable property of \( \mathcal{P}^{1}_{N,\epsilon} \) that does not appear in \( \mathcal{P}^{2}_{N,\epsilon} \). This characterization provides closed-form solutions to a distributionally robust newsvendor problem with a single scenario aunder both \( \mathcal{P}^{1}_{N,\epsilon} \) and \( \mathcal{P}^{2}_{N,\epsilon} \), suggesting that \( \tilde{x}^{2}_{N}(\epsilon) \) is better informed than \( \tilde{x}^{1}_{N}(\epsilon) \); specifically, unlike \( \tilde{x}^{1}_{N}(\epsilon) \), \( \tilde{x}^{2}_{N}(\epsilon) \) accounts for both the critical fractile and the Wasserstein radius \( \epsilon \). Through a penalty-based interpretation of the dual worst-case expectation problem, we further reveal the superior out-of-sample performance of 2-Wasserstein balls for general \( \Xi \), which offer a wider range of radii yielding solutions that outperform SAA and RO, aligning with the observation in Figure \ref{fig:oos}.

\subsection{Notation}
\label{sec:notations}
Throughout this paper, let $q$ be a scalar such that $1/p+1/q=1$, ensuring $\|\cdot\|_q$ is the dual norm of $\|\cdot\|$. For an integer $n$, define $[n]:=\{1,\cdots, n\}$ and let $\mathbb R^n_+$ represent the $n$-dimensional nonnegative orthant. For a set $\mc K \subseteq \mathbb R^n$, $\mc K^*$ denotes its dual cone, that is, $\mc K^* = \{y \in \mb R^n: y^T x \ge 0, \ \forall x \in \mc K\}$, {\color{black}and $\mbox{recc}(\mc K)$ represents the recession cone of $\mc K$, that is $\mbox{recc}(\mc K)=\{\vec u \in \mathbb R^n: \forall u \in \mathcal K, u+ \alpha \vec u \in \mc K, \ \forall \alpha \in \mathbb R_+\}$}. For a proper cone $\mc K$, $\preceq_{\mc K}$ denotes the generalized inequality defined by $\mc K$, i.e., $x \preceq_{\mc K} y \Leftrightarrow y-x \in \mc K$. For a vector $v \in \mathbb{R}^k$ and a scalar $a$, $v^a$ denotes the componentwise exponentiation of $v$ with exponent $a$, i.e., $v^a = (v_j^a)_{j \in [k]}$, and $|v|$ denotes the componentwise absolutization of $v$, that is, $|v| = (|v_j|)_{j \in [k]}$.


\section{Preliminaries and assumptions}
\noindent The dual of \eqref{prob:second-stage} is given as
    \(
        \max_{\pi \ge 0} \left\{(Tx + \xi)^\top \pi : W^\top \pi \preceq_{\mc K_y^*} q\right\}.
    \)
  We let $\Pi$ denote the dual feasible region of \eqref{prob:second-stage}; that  is, $\Pi = \{\pi \ge 0:W^\top \pi \preceq_{\mc K_y^*} q\}$. 
\begin{assu}
\begin{enumerate}
    \item[(i)] Complete recourse: \eqref{prob:second-stage} is feasible for any right-hand side $Tx + \xi$ of the affine constraints. 
    \item[(i)] Recourse with dual strict feasibility: $\Pi$ has an interior point.
  \end{enumerate}\label{assum}
\end{assu}
\noindent Assumption \ref{assum}(i) ensures $\Pi$ is bounded, which can be achieved approximately by adding slack variables with large penalties. 
Additionally, $\Pi$ is closed since it is the intersection of a nonnegative orthant and the inverse image of a closed set $\mathcal{K}^*_y$ under an affine mapping, and thus $\Pi$ compact. Assumption \ref{assum}(ii) ensures $Z(x, \xi) > -\infty,$ $\forall x \in \mathbb{R}^{n_x},\xi \in \Xi$, guaranteeing strong duality of the second-stage problem. 

We make the following additional common assumption on $\Xi$:
\begin{assu}
    The support set $\Xi$ is nonempty, closed, and convex.
  \label{assum:support}
\end{assu}
\noindent The following proposition is well-established in the literature \cite{hanasusanto2018conic,blanchet2019quantifying,gao2023distributionally}:
\begin{prop} Under Assumption \ref{assum}, for any fixed $x \in \mb R^{n_x}$, 
  the dual of $\mathbb E_{\mathbb P \in \mathcal P^{r}_{N,\epsilon}}[Z(x,\tilde \xi)]$ is
\begin{equation} \min_{\lambda \ge 0} {\epsilon^r} \lambda + \frac{1}{N} \sum_{i\in [N]} \sup_{\xi \in \Xi} Z(x,\xi) - \lambda \|\xi-\zeta^i\|^r,\label{prob:worst-case-expectation-dual}\end{equation}
and when ${\epsilon} > 0$, strong duality holds and the minimum is attained.
\label{prop:DRO-duality}
\end{prop}

Let $\mathcal D$ denote the set of all normalized feasible rays of $\Xi$ with respect to the $l_p$-norm, that is $\mathcal D :=\{d \in \mathbb R^k: d \in \mbox{recc}\Xi, \|d\| = 1\}$. It is observed in \cite{byeon2022two} (see Corollary 3.4. of the paper), 
\begin{equation}\sup_{\xi \in \Xi} Z(x,\xi) - \lambda \|\xi-\zeta^i\| < \infty,  \ \forall i \in [N] \Longleftrightarrow \lambda \ge \max_{\pi \in \Pi, d \in \mathcal D} \pi^\top d,\label{eq:r1:lambda-bound}\end{equation}
imposing an implicit lower bound on $\lambda$ when $r=1$, which is well-defined since $\Pi \times \mathcal D$ is compact. Let \( L := \max_{\pi \in \Pi} \|\pi\|_q \) where $q$ is such that $1/p+1/q=1$. It is also shown in \cite{byeon2022two} (see Lemma 2.2 of the paper) that $Z(x,\cdot)$ is $L$-Lipschitz continuous with respect to the $l_p$-norm for any $x \in \mathbb R^{n_x}$ under Assumption \ref{assum}.

\section{$\Xi=\mathbb R^k$ or $\Xi=\mathbb R^k_+$}
We first analyze the cases where $\Xi = \mathbb{R}^k$ or $\mathbb{R}^k_+$ and show that $\mathcal{P}^1_{N,\epsilon}$ contains a pathological worst-case distribution. Consequently, $\tilde{x}^1_N(\epsilon)$ mirrors $\tilde{x}^\texttt{SAA}_N$, regardless of $\epsilon$. In contrast, TSDRO over $\mathcal{P}^2_{N,\epsilon}$ considers a worst-case distribution that redistributes probability masses at sample points in a regularized manner, guided by the optimal Lagrange multiplier $\lambda^*$, which varies with $\epsilon$. This highlights the superiority of $x^2_N(\epsilon)$ over $x^1_N(\epsilon)$, as demonstrated through a newsvendor example.

\subsection{$r=1$}
Note first that when \( \Xi = \mathbb{R}^k \), the set \( \mathcal{D} \) in \eqref{eq:r1:lambda-bound} is given by \( \mathcal{D} = \{d \in \mathbb{R}^k : \|d\| = 1\} \). In this case, the implicit lower bound on $\lambda$ is 
\(
\max_{\pi \in \Pi, d \in \mathcal{D}} \pi^\top d = \max_{\pi \in \Pi} \|\pi\|_q =: L,
\)  
since the dual norm of the \( l_p \)-norm is the \( l_q \)-norm. When \( \Xi = \mathbb{R}^k_+ \), i.e., \( \mathcal{D} = \{d \in \mathbb{R}^k_+ : \|d\| = 1\} \), we similarly have  
\(
\max_{\pi \in \Pi, d \in \mathcal{D}} \pi^\top d = L.
\)  
This is because, for any \( \pi \in \Pi \setminus \{0\} \subseteq \mathbb{R}^k_+ \), we have $\max_{d \in \mathcal D}\pi^\top d = \|\pi\|_q$ since \( d = \frac{\pi^{q-1}}{\|\pi\|_q^{q-1}}  \in \mathcal{D} \), where \( \pi^{q-1} \) denotes the componentwise exponentiation of \( \pi \) with exponent \( q-1 \).
 This gives the following proposition with its proof provided in Appendix \ref{appendix:1-wass:sol}:
\begin{prop}\label{prop:r1:worst-expectation}
When \( \Xi = \mathbb{R}^k \) or \( \Xi = \mathbb{R}^k_+ \), for any \( \epsilon > 0 \),
\[
\sup_{\mathbb{P} \in \mathcal{P}^{1}_{N,\epsilon}} \mathbb{E}_{\mathbb{P}}[Z(x, \tilde{\xi})] = \epsilon L + \frac{1}{N} \sum_{i \in [N]} Z(x, \zeta^i).
\]
\end{prop}
\noindent Proposition \ref{prop:r1:worst-expectation} implies that 
\(\tilde{x}^{1}_{N}(\epsilon)\) is \emph{\(\epsilon\)-invariant and always coincides with \(\tilde{x}^{\texttt{SAA}}_N\)} unless \(\epsilon = \infty\).  
This occurs because \(\mathcal{P}^{1}_{N,\epsilon}\) includes a pathological probability distribution that shifts a diminishing fraction of probability mass away from a single sample point in the direction of
\begin{equation}
d_{\pi^*} := \frac{(\pi^*)^{q-1}}{\|\pi^*\|_q^{q-1}} \label{r1:direction}
\end{equation}
for some \(\pi^* \in \Pi^* := \arg\max_{\pi \in \Pi} \|\pi\|_q\). Note that \( \|d_{\pi^*}\| = 1 \). This is formally stated in the following proposition, the proof of which is provided in Appendix \ref{appendix:proof:1-wass-worst-dist}:
\begin{prop}\label{prop:1-wass-worst-dist}
Let $\Xi = \mathbb R^k$ or $\mathbb R^k_{+}$. 
For any $\epsilon >0$, the worst-case expectation $\sup_{\mathbb P \in \mathcal P^{1}_{N,\epsilon}}\mathbb E_{\mathbb P}[Z(x,\tilde \xi)]$ is either attained or asymptotically attained by a probability distribution that reallocates probability mass on a single sample point $\zeta^{i'}$ for some $i' \in [N]$ in the direction of $d_{\pi^*}$, as defined in \eqref{r1:direction}. Define $\Pi^*_{i,x} := \arg\max_{\pi \in \Pi} \pi^\top(Tx + \zeta^i)$.
\begin{enumerate}
\item[(i)] \textbf{Attainable case:} If $\exists i' \in [N]: \Pi^* \cap \Pi^*_{i',x} \neq \emptyset$, then for any $\pi^* \in \Pi^* \cap \Pi^*_{i',x}$, $\sup_{\mathbb P \in \mathcal P^{1}_{N,\epsilon}}\mathbb E_{\mathbb P}[Z(x,\tilde \xi)]$ is attained by the following distribution in $\mathcal P^{1}_{N,\epsilon}$:
\begin{equation}\mathbb Q_x = \frac{1}{N}\left(\sum_{i \in [N]: i \neq i'} \delta_{\zeta^i} + \delta_{\zeta^{i'}+ \epsilon Nd_{\pi^*}}\right), \label{eq:dist:r1:R}\end{equation}
which moves the $\frac{1}{N}$-mass at $\zeta^{i'}$ to $\zeta^{i'}+{\epsilon N}d_{\pi^*}$.
\item[(ii)] \textbf{Asymptotically attainable case:} Otherwise, if $\Pi^* \cap \Pi^*_{i,x} = \emptyset$ for all $i \in [N]$, $\sup_{\mathbb P \in \mathcal P^{1}_{N,\epsilon}}\mathbb E_{\mathbb P}[Z(x,\tilde \xi)]$ can only be achieved asymptotically by a sequence of distributions $\{\mathbb Q_{x,\Delta}\}_{\Delta} \subseteq \mathcal P^{1}_{N,\epsilon}$ as $\Delta \rightarrow 0$, given by: 
 \begin{equation}\mathbb Q_{x,\Delta} = 
\frac{1}{N}\sum_{i \in [N]: i \neq i'} \delta_{\zeta^i} + \Delta\delta_{\zeta^{i'}+ \frac{\epsilon}{\Delta} d_{\pi^*}} + \left(\frac{1}{N} - \Delta\right)\delta_{\zeta^{i'}},\label{eq:dist:r1:R:2}\end{equation}
for some $i' \in [N]$ and $\pi^* \in \Pi^*$. Here, $\Delta \in (0,\frac{1}{N}]$ represents the probability mass transported from $\zeta^{i'}$ in the direction of $d_{\pi^*}$ 
by the $\frac{\epsilon}{\Delta}$-distance. 
\end{enumerate}
\end{prop}
\begin{rema}
Proposition \ref{prop:r1:worst-expectation} was previously presented in \cite{hanasusanto2018conic,esfahani2018data,byeon2022two} only for \( \Xi = \mathbb{R}^k \). In \cite{duque2022distributionally}, it was claimed that a similar argument holds for a convex conic \( \Xi \); however, a counterexample to this claim was provided in \cite{byeon2022two}. In addition, to the best of the author's knowledge, Proposition \ref{prop:1-wass-worst-dist} is the first to precisely characterize the worst-case distribution for TSDRO with right-hand side uncertainty over 1-Wasserstein balls for both unconstrained and nonnegative orthant supports. It also aligns with the result in \cite{gao2023distributionally}, which states that if the worst-case distribution exists, it supports at most \( N+1 \) points. Based on the proofs of Propositions \ref{prop:r1:worst-expectation} and \ref{prop:1-wass-worst-dist}, it can be seen that both propositions also hold for other conic supports \( \Xi \), provided that \( (\pi^*)^{q-1} \in \mathrm{recc}(\Xi) \) for some \( \pi^* \in \Pi^* \), thereby ensuring \( d_{\pi^*} \in \mathcal{D} \) and \( \max_{\pi \in \Pi, d \in \mathcal{D}} \pi^\top d = L \).
\end{rema}

\subsection{$r=2$}
Conversely, when $r=2$, we demonstrate that the worst-case probability distribution allocates probability masses in a regularized manner, governed by the optimal Lagrange multiplier $\lambda^*$, associated with the Wasserstein constraint. It is formally stated in the following proposition, the proof of which is in Appendix \ref{appendix:proof:2-wass-worst-dist}:
\begin{prop}\label{prop:2-wass-worst-dist}
When $\Xi = \mathbb R^k$ or $\mathbb R_+^k$, the optimal Lagrange multiplier $\lambda^*$ of $\sup_{\mathbb P \in \mathcal P^{2}_{N,\epsilon}}\mathbb E_{\mathbb P}[Z(x,\tilde \xi)]$ should satisfy
\begin{equation}
    \frac{1}{N} \sum_{i \in [N]} \sum_{\pi \in \Pi^*_{i, x, \lambda^*}} \mu^i_{\pi} \frac{\|\pi\|_q^2}{4(\lambda^*)^2} = \epsilon^2 \label{eq:worst-case-prob-dist:r2:dual}
\end{equation}
for some $(\mu^i_{\pi})_{\pi \in \Pi^*_{i,x,\lambda^*}} \ge 0$ satisfying $\sum_{\pi \in \Pi^*_{i,x,\lambda^*}} \mu^i_{\pi}=1, \forall i \in [N]$, 
where $\Pi^*_{i,x,\lambda^*}:= \arg\max_{\pi \in \Pi} \pi^\top (Tx + \zeta^i) + \frac{1}{4\lambda^*}\|\pi\|_q^2$. In addition, $\sup_{\mathbb P \in \mathcal P^{2}_{N,\epsilon}}\mathbb E_{\mathbb P}[Z(x,\tilde \xi)]$ is attained by a probability distribution 
\begin{equation}\mathbb Q_x = \frac{1}{N}\sum_{i \in [N]}\sum_{\pi \in \Pi^*_{i,x,\lambda^*}}\mu^i_{\pi}\delta_{\zeta^i+\frac{\|\pi\|_q}{2\lambda^*}{d_\pi}}.\label{eq:worst-prob-dist:r2}\end{equation}
\end{prop}

\begin{figure}\centering
\begin{subfigure}[b]{0.33\textwidth}
      \centering      
      \includegraphics[width=\textwidth]{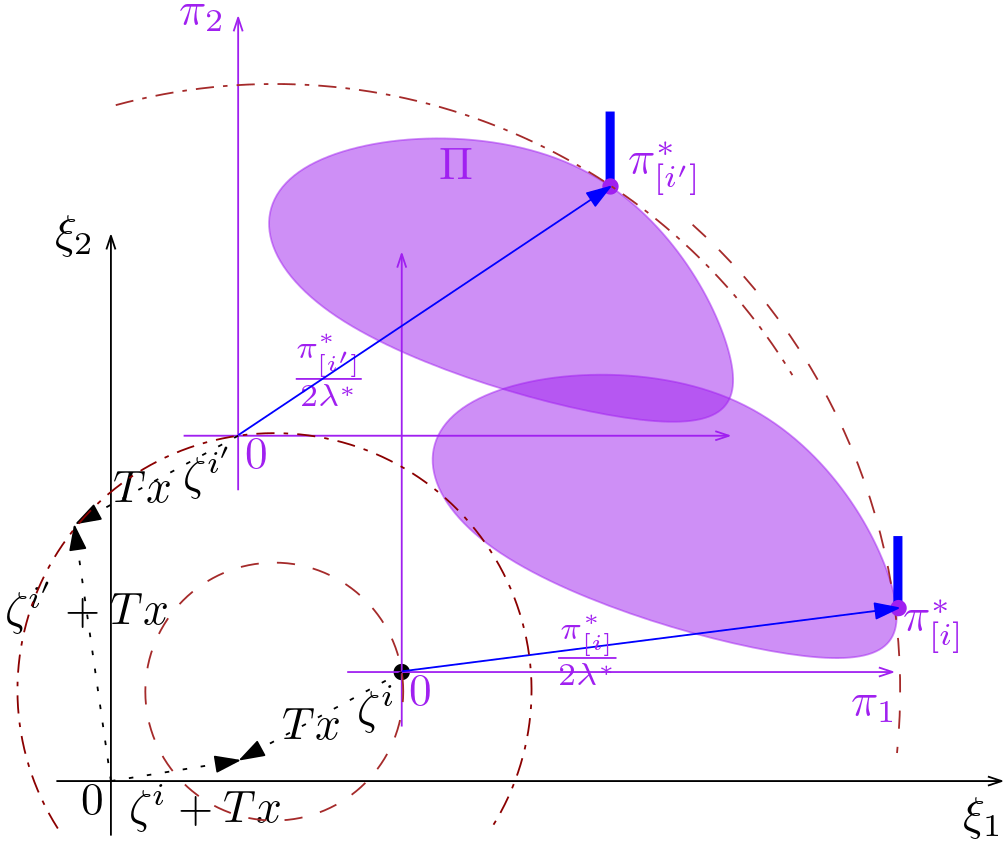}
      \caption{$\lambda^*=0.5$}\label{fig:worst-case-dist-r2:1}
  \end{subfigure} 
    \hspace{1cm}
  \begin{subfigure}[b]{0.33\textwidth}
      \centering      
      \includegraphics[width=\textwidth]{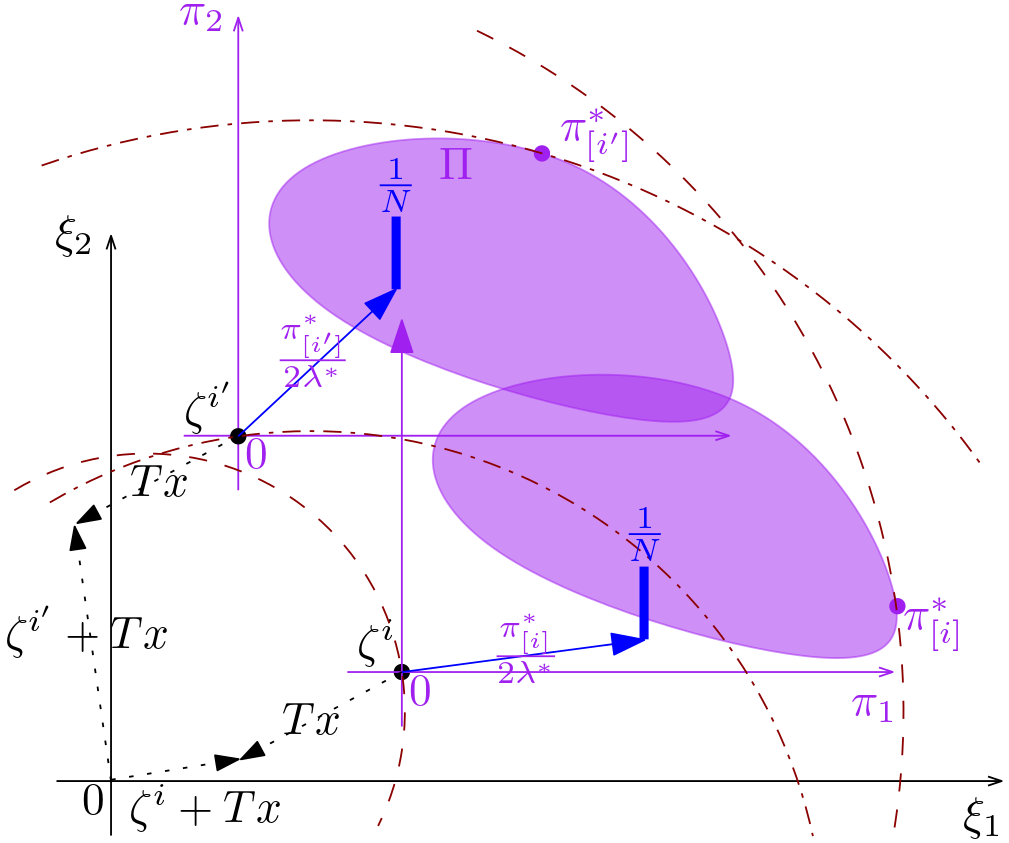}
      \caption{$\lambda^*=1$}\label{fig:worst-case-dist-r2:2}
  \end{subfigure} 
  \caption{Worst-case probability distribution in $\mathcal P^{2}_{N,\epsilon}$, where $\pi^*_{[i]} = \arg\max_{\pi \in \Pi}\pi^\top (Tx+\zeta^i)+\frac{1}{4\lambda^*}\|\pi\|^2$.} 
  \label{fig:worst-case-dist-r2}
\end{figure}

Figure \ref{fig:worst-case-dist-r2} illustrates the worst-case probability distribution for \( r = 2,p=2, \) and \( \Xi = \mathbb{R}^k \) when $\lambda^*=0.5$ or 1. In each figure, the \(\pi\)-space is overlaid onto the \(\xi\)-space, with each origin positioned at \(\zeta^i\) for \( i \in [N] \). Note that \( o_i(\pi) := \pi^\top(Tx + \zeta^i) + \frac{1}{4\lambda^*}\|\pi\|_2^2 \) is a strictly convex quadratic function with circular contours centered at \( -2\lambda^*(Tx + \zeta^i) \). The dashed-dot circle (and the dashed circle) represents the contour of \( o_{i'}(\pi) \) (and \( o_i(\pi) \)). The point \( \pi_{[i]}^* \) denotes the \(\pi\) that maximizes \( o_i(\pi) \) over \(\Pi\). 
Since each \( \max_{\pi \in \Pi} o_i(\pi) \) has a unique optimum, the worst-case distribution shifts the entire probability mass at \( \zeta^i \) to \( \zeta^i + \frac{1}{2\lambda^*} \pi_{[i]}^* \). 
If there were multiple optimizers, the \( \frac{1}{N} \)-probability mass could be distributed across them, each scaled by \( 2\lambda^* \). 
If \( N = 2 \), the first case corresponds to the case where 
\(
\epsilon = \sqrt{1/2} \| ( \pi^*_{[1]}; \pi^*_{[2]} ) \|
\)
and the second case to 
\(
\epsilon = \sqrt{1/8} \| ( \pi^*_{[1]}; \pi^*_{[2]} ) \|,
\)
with the associated \( \pi^* \). The $\epsilon$-dependent optimal \( \lambda \) value, as in \eqref{eq:worst-case-prob-dist:r2:dual}, determines the worst-case distribution; (i) \( \lambda^* \) adjusts the center of each circular contour, thereby changing \( \pi^*_{[i]} \); (ii) a smaller \( \lambda^* \) results in a greater shift distance.

\begin{rema}
Proposition \ref{prop:2-wass-worst-dist} implies existence of the worst-case distribution in $\mathcal P^2_{N,\epsilon}$. Generally, the existence is not guaranteed for unbounded $\Xi$, and several conditions ensuring the existence were proposed in \cite{blanchet2019quantifying,gao2023distributionally,esfahani2018data,yue2022linear}. One can see that $\sup_{\mathcal P^2_{N,\epsilon}}\mathbb E_{\mathbb P}[Z(x,\tilde \xi)]$ indeed
satisfies the existence condition in \cite{yue2022linear} for any $x$ under Assumption \ref{assum}.
\end{rema}


\subsection{Example: single-scenario newsvendor problem}
Consider a classic newsvendor problem. Let \( \texttt{c} \) denote the unit cost per newspaper and \( \texttt{p} \) the selling price per newspaper, where \( 0 < \texttt{c} < \texttt{p} \). Define \( \tilde{\xi} \) as the negative of the uncertain demand, for which a single sample \( \zeta \) is observed. Therefore, we omit the sample index \( i \) in this section. The here-and-now decision \( x \) represents the number of newspapers to order, while the wait-and-see decision \( y \) denotes the quantity of newspapers leftover. The distributionally robust newsvendor problem over $\mathcal P^{r}_{N,\epsilon}$ can be formulated as follows:
\begin{equation}
\min_{x \geq 0} \ (\texttt{c} - \texttt{p}) x + \sup_{\mathbb{P} \in \mathcal{P}^{r}_{N,\epsilon}} \mathbb{E}_{\mathbb{P}} [Z(x, \tilde{\xi})],
\label{prob:newsvendor}\end{equation}
where
\(
Z(x, \tilde{\xi}) = \texttt{p} \max\{x + \tilde{\xi}, 0\}= \min \{ \texttt{p} y : y \geq x + \tilde{\xi}, \, y \geq 0 \}
\)
represents unfulfilled sales if the demand, \( -\tilde{\xi} \), is lower than the order quantity \( x \). The optimal solution to \eqref{prob:newsvendor} can be characterized as follows, the proof of which is in Appendix \ref{appendix:proof:newsvendor}:
\begin{prop}
Let $\Xi = \mathbb{R}$, $N=1$, and let $\zeta$ denote the sample point. For any $\epsilon \in (0,\infty)$, the optimal solution to \eqref{prob:newsvendor} when $r=1$ and $2$ is, respectively:  
\[
{x}^{1}_{1}(\epsilon) = -\zeta, \quad 
{x}^{2}_{1}(\epsilon) = \max\left\{0, -\zeta - \frac{\epsilon}{2} \sqrt{\frac{\texttt{p}}{\texttt{p} - \texttt{c}}}\right\}.
\]\label{prop:newsvendor}
\end{prop}
\noindent Proposition \ref{prop:newsvendor} considers the case where \( \tilde{\xi} \) is supported on \( \Xi = \mathbb{R} \), such as when \( \tilde{\xi} \) follows a Gaussian distribution. This implies that negative demand can occur, meaning the newsvendor may purchase newspapers brought in by customers at the selling price \( \texttt{p} \). Hence, the newsvendor must account for this potential supply from customers when making the here-and-now decision. 

\begin{rema}
Proposition \ref{prop:newsvendor} shows that ${x}^2_{1}(\epsilon)$ depends on the sample point \(\zeta\), the risk appetite \(\epsilon\), and the critical fractile \(\frac{\texttt{p} - \texttt{c}}{\texttt{p}}\). Given a true cumulative distribution function $F$ of $-\tilde \xi$, the optimal quantity \(x^*\) is known to be \(F^{-1}\left(\frac{\texttt{p} - \texttt{c}}{\texttt{p}}\right)\).
For \(r=2\), a lower critical fractile results in purchasing fewer newspapers, similar to \(x^*\), for \(\epsilon \in \left(0, -2\zeta \sqrt{\frac{\texttt{p} - \texttt{c}}{\texttt{p}}}\right)\). As \(\epsilon\) increases, ${x}^2_{1}(\epsilon)$ decreases and converges to the RO solution \({\mathcal X}_{\texttt{RO}} = \{0\}\), as \(\sup_{\xi \in \mathbb{R}} \texttt{p} \max\{x + \xi, 0\} = \infty\) for any \(x\). Conversely, ${x}^1_{1}(\epsilon)$ remains independent of \(\epsilon\) and the critical fractile, highlighting that \(\mathcal{P}^2_{N,\epsilon}\) enables more informed decisions compared to \(\mathcal{P}^1_{N,\epsilon}\). The out-of-sample performance of a decision \(x^r_1(\epsilon)\) is given by:  
\(
O^r_1(\epsilon) = (\texttt{c} - \texttt{p}) x^r_1(\epsilon) + \int_{\xi \in \mathbb{R}} \texttt{p} \max\{x^r_1(\epsilon) + \xi, 0\} d\mathbb{P}^{\texttt{true}}(\xi),
\)
which simplifies to:  
\(
O^r_1(\epsilon) = (\texttt{c} - \texttt{p}) x^r_1(\epsilon) + \texttt{p} x^r_1(\epsilon) \mathbb{P}^{\texttt{true}}[x^r_1(\epsilon) + \xi \geq 0] + \texttt{p} \int_{-x^r_1(\epsilon)}^\infty \xi d\mathbb{P}^{\texttt{true}}(\xi).
\)

Figure~\ref{fig:newsvendor} shows the 10th and 90th percentiles, as well as the average, of \(\tilde{O}^r_1(\epsilon)\) and \(\tilde{x}^2_1(\epsilon)\) for a specific instance over 2500 simulation runs for each \(\epsilon\) values ranging from 0.005 to 15. The results indicate an \(\epsilon\)-value where the average \(\tilde{x}^2_1(\epsilon) \approx x^*\), achieving consistently high revenues even with a small sample size (\(N=1\)) and low-margin setting where \(\texttt{p}/\texttt{c}\) is close to 1, in contrast to \(\tilde{x}^1_{1}(\epsilon)\). Notably, $\tilde{x}^1(\epsilon) = \tilde{x}^2(0)$ for all $\epsilon > 0$, and it fails to include $x^*$ within the 80\% confidence band.

\begin{figure}[h]
    \centering              
    \begin{subfigure}[b]{0.33\textwidth}
        \centering
        \includegraphics[width=\textwidth]{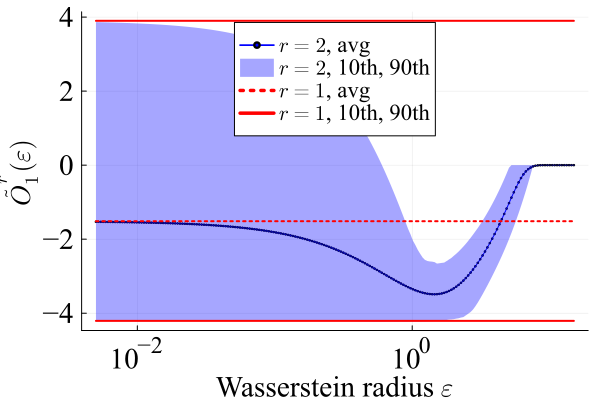}
    \end{subfigure}
    \begin{subfigure}[b]{0.33\textwidth}
        \centering
        \includegraphics[width=\textwidth]{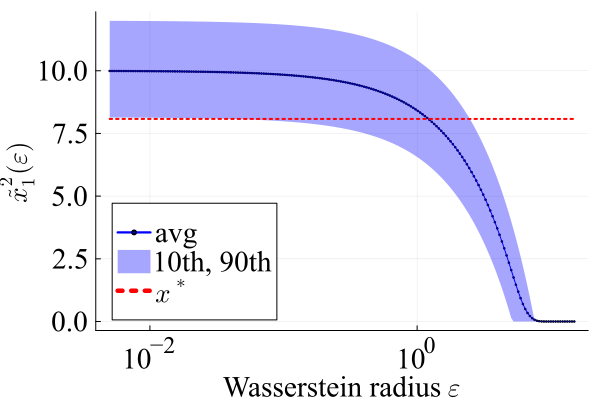}
    \end{subfigure}
    \caption{\(\tilde{O}^r_{1}(\epsilon)\) for \(r=1,2\), and \(\tilde{x}^2_1(\epsilon)\) when \(\mathbb{P}^{\texttt{true}} = \mathcal{N}(-10, 1.5)\), \(\texttt{p} = 5\), \(\texttt{c} = 4.5\).}
    \label{fig:newsvendor}
\end{figure}
\end{rema}


\section{Discussion on general $\Xi$}
\noindent Built upon Proposition \ref{prop:DRO-duality}, for any \(\epsilon > 0\), \eqref{prob:TSDRO} can be formulated as  
\(
\min_{\lambda \geq 0} \epsilon^r \lambda + \min_{x \in \mathcal{X}} c^\top x + \frac{1}{N} \sum_{i \in [N]} \sup_{\xi \in \Xi} \big(Z(x, \xi) - \lambda \|\xi - \zeta^i\|^r\big),
\)  
where the minimization over \(\lambda\) and \(x\) is interchanged. For a given \(\lambda > 0\), the inner minimization seeks optimal \(x\) with respect to the modified risk measure:
\(
\frac{1}{N} \sum_{i \in [N]} g_{i}^{r}(x, \lambda),
\) where
\(
g_{i}^{r}(x, \lambda) := \sup_{\xi \in \Xi} Z(x, \xi) - \lambda \|\xi - \zeta^i\|^r.  
\)
Note $g_{i}^{r}(x, \lambda)$ can be interpreted as a penalized version of \(
Z(x, \zeta^i) = \sup_{\xi \in \Xi} \big\{ Z(x, \xi) : \xi = \zeta^i \big\},
\)
where a penalty is imposed for violating the constraint \(\xi = \zeta^i\) with a penalty parameter \(\lambda > 0\). When \(r = 1\), this corresponds to an \emph{exact} penalty, whereas for \(r = 2\), it represents an \emph{inexact} quadratic penalty.  

It is well-known in the literature that for exact penalty methods, a sufficiently large penalty parameter ensures equivalence between the penalized and non-penalized problems. Specifically, for an \(L\)-Lipschitz continuous function \(Z(x, \cdot)\), any \(\lambda \geq L\) guarantees equivalence, i.e., \(g_{i}^{1}(x, \lambda) = Z(x, \zeta^i)\) for any $\lambda \ge L$. For quadratic penalties, equivalence requires \(\lambda \to \infty\), meaning \(g_{i}^{2}(x, \lambda) \neq Z(x, \zeta^i)\) unless \(\lambda \to \infty\). Consequently, when \(r = 1\), for a certain range of small values of \(\epsilon\), the optimal \(\lambda\) can be large enough that the solution coincides with \(\tilde x^{\texttt{SAA}}_N\). In contrast, \emph{when \(r = 2\), the solution can deviate from \(\tilde x^{\texttt{SAA}}_N\) as soon as \(\epsilon > 0\).}

Furthermore, when \(\Xi\) is compact such that \(\sup_{\xi \in \Xi} Z(x, \xi) < \infty\) for all \(x \in \mathcal{X}\), the supremum \(\sup_{\xi \in \Xi} Z(x, \xi)\) is attained at a boundary point \(\xi^{\texttt{ro}}\) of \(\Xi\). For large $\epsilon$ values such that \(\mathcal{P}^{r}_{N,\epsilon}\) includes \(\delta_{\xi^{\texttt{ro}}}\), the probability distribution that concentrates all probability mass on \(\xi^{\texttt{ro}}\), the decision reduces to a solution \(x^{\texttt{RO}} \in {\mathcal X}_{\texttt{RO}}\), i.e.,  
\(
\tilde x^{r}_N(\epsilon) = x^{\texttt{RO}} \in \arg\min_{x \in \mathcal{X}} c^\top x + \sup_{\xi \in \Xi} Z(x, \xi).
\)  
Notably, \(\mathcal{P}^{1}_{N,\epsilon}\) includes \(\delta_{\xi^{\texttt{ro}}}\) for smaller values of \(\epsilon\) compared to \(\mathcal{P}^{2}_{N,\epsilon}\). This is because, for \(a = [\|\xi^{\texttt{ro}} - \zeta^i\|]_{i \in [N]}\), we have  
\(
W^{2}(\delta_{\xi^{\texttt{ro}}}, \empi) = \sqrt{\frac{1}{N} \sum_{i \in [N]} \|\xi^{\texttt{ro}} - \zeta^i\|^2} = \sqrt{\frac{1}{N}} \|a\|_2 \geq \frac{1}{N} \|a\|_1 = W^{1}(\delta_{\xi^{\texttt{ro}}}, \empi),
\)
where equality holds only when all elements of \(a\) are equal.  

This implies that  
\(
\{\epsilon : \tilde x^{1}_N(\epsilon) \neq \tilde x^{\texttt{SAA}}_N \text{ and } \tilde x^{1}_N(\epsilon) \not\in {\mathcal X}_{\texttt{RO}} \mbox{ a.s.}\} \subseteq \{\epsilon : \tilde x^{2}_N(\epsilon) \neq \tilde x^{\texttt{SAA}}_N \text{ and } \tilde x^{2}_N(\epsilon) \not\in {\mathcal X}_{\texttt{RO}}\mbox{ a.s.}\},
\)  
and the numerical result on facility location problem instances in Figure \ref{fig:oos} empirically demonstrate this implication. Given that finding the optimal \(\epsilon\) value is challenging, having a wider range of \(\epsilon\) values for which the DRO solution outperforms the SAA and RO counterparts highlights the advantage of TSDRO with $\mathcal P^{2}_{N,\epsilon}$ over its $\mathcal P^{1}_{N,\epsilon}$ counterpart.

\appendix
\vspace{-4mm}
\section{Proof of Proposition \ref{prop:r1:worst-expectation}}\label{appendix:1-wass:sol}
From \eqref{eq:r1:lambda-bound}, we have $\lambda \ge L$, and from the $L$-Lipschitz continuity of $Z(x,\cdot)$, we have
\(
\sup_{\xi \in \Xi} Z(x, \xi) - \lambda \|\xi - \zeta^i\| = Z(x, \zeta^i) + \sup_{\xi \in \Xi} Z(x, \xi) - Z(x, \zeta^i) - \lambda \|\xi - \zeta^i\| \le Z(x, \zeta^i) + \sup_{\xi \in \Xi}  (L- \lambda) \|\xi - \zeta^i\| =  Z(x, \zeta^i).
\)
Therefore, the problem becomes $\min_{\lambda \ge L} \epsilon \lambda + \frac{1}{N}\sum_{i \in [N]}Z(x,\zeta^i) = \epsilon L + \frac{1}{N}\sum_{i \in [N]}Z(x,\zeta^i)$, so the desired result follows.
\qed

\section{Proof of Proposition \ref{prop:1-wass-worst-dist}}
\label{appendix:proof:1-wass-worst-dist}
Let us first consider the case where $\exists i' \in [N]: \Pi^* \cap \Pi^*_{i',x} \neq \emptyset$. 
It is easy to verify that the probability distribution $\mathbb Q_x$ defined as \eqref{eq:dist:r1:R} 
is in $\mathcal P^{1}_{N,\epsilon}$, since
$\|d_{\pi^*}\| = 1$, and it evaluates the objective function as
\begin{subequations}\tiny
\begin{align}
\mathbb E_{\mathbb Q_x}[Z(x,\tilde \xi)] & = \frac{1}{N} \left(\sum_{i \in [N], i \neq i'}Z(x,\zeta^i) + Z(x,\zeta^{i'}+ \epsilon Nd_{\pi^*})\right) \\
&= \frac{1}{N} \sum_{i \in [N]}Z(x,\zeta^i) + \frac{1}{N}\left(Z(x,\zeta^{i'}+ \epsilon Nd_{\pi^*}) - Z(x,\zeta^{i'})\right)\\
& \ge \frac{1}{N} \sum_{i \in [N]}Z(x,\zeta^i) + \frac{1}{N}\left((\pi^*)^\top (Tx + \zeta^{i'}+ \epsilon Nd_{\pi^*}) - Z(x,\zeta^{i'})\right) \label{eq:worst-q:r1:1}\\
& = \frac{1}{N} \sum_{i \in [N]}Z(x,\zeta^i) + \frac{1}{N}\left((\pi^*)^\top \epsilon Nd_{\pi^*}\right) \label{eq:worst-q:r1:2}\\
& = \frac{1}{N} \sum_{i \in [N]}Z(x,\zeta^i) + \frac{1}{N}\left(\epsilon N\frac{\|\pi^*\|_q^q}{\|\pi^*\|_q^{q-1}}\right) = \frac{1}{N} \sum_{i \in [N]}Z(x,\zeta^i) + \epsilon L\label{eq:worst-q:r1:3}
\end{align}
\end{subequations}
where \eqref{eq:worst-q:r1:1} follows from $\pi^* \in \Pi$ and $Z(x,\zeta^{i'}+ \epsilon Nd_{\pi^*})=\max_{\pi \in \Pi} \pi^\top (Tx + \zeta^{i'}+ \epsilon Nd_{\pi^*})$, \eqref{eq:worst-q:r1:2} from $\pi^* \in \Pi^*_{i',x}$, \eqref{eq:worst-q:r1:3} from the definition of $d_{\pi^*}$ and from $\pi^* \in \Pi^*$. Therefore, $\mathbb Q_x$ attains the worst-case expectation. 

For Case (ii), note that for any $i \in [N]$ and $\xi \in \Xi\setminus \{\zeta^i\}$, we have
\begin{subequations}\tiny
\begin{align}
{Z(x,\xi) - Z(x,\zeta^i)} 
  &=  \max_{\pi \in \Pi }\pi^\top (Tx + \zeta^i + \xi - \zeta^i) - Z(x,\zeta^i)\\
& \le \max_{\pi \in \Pi }\pi^\top (Tx + \zeta^i) + \max_{\pi \in \Pi }\pi^\top |\xi - \zeta^i| -  Z(x,\zeta^i) = {\max_{\pi \in \Pi }\pi^\top |\xi - \zeta^i|}\label{eq:worst-q:r1:5}\\
& \le \max_{\pi \in \Pi}\|\pi\|_q \|\xi - \zeta^i\| = L \|\xi - \zeta^i\|\label{eq:worst-q:r1:6},
\end{align}    
\end{subequations}
where the inequality in \eqref{eq:worst-q:r1:5} follow from $\max_{a \in \mathcal A}f(a) + g(a) \le \max_{a \in \mathcal A}f(a) + \max_{a \in \mathcal A} g(a)$ and $\Pi \subseteq \mathbb R^k_+$ and the equality in \eqref{eq:worst-q:r1:5} is from $\max_{\pi \in \Pi }\pi^\top (Tx + \zeta^i) = Z(x,\zeta^i)$; \eqref{eq:worst-q:r1:6} from Hölder's inequality. 
Note that if $\exists \pi^* \in \Pi^*$ and a scalar $\texttt{c}$ such that $|\xi_j - \zeta^i_j| = \texttt{c} (\pi^*_j)^{q-1}$ for all $j \in [k]$, then  
\eqref{eq:worst-q:r1:5} should hold with strict inequality since $\max_{\pi \in \Pi}\pi^\top |\xi - \zeta^i| = \texttt{c}\max_{\pi \in \Pi}\pi^\top (\pi^*)^{q-1} \le \texttt{c}\max_{\pi \in \Pi}\|\pi\|_q \|(\pi^*)^{q-1}\| = \texttt{c}\max_{\pi \in \Pi}\|\pi\|_q \|\pi^*\|^{q-1}_q = \texttt{c}L$ is achieved by $\pi^*$ and $\Pi^* \cap \Pi^*_{i,x} =\emptyset$. Otherwise, \eqref{eq:worst-q:r1:6} holds with strict inequality. Therefore, we have for any $\xi \in \Xi \setminus\{\zeta^i\}$,
\begin{equation}\tiny{Z(x,\xi) - Z(x,\zeta^i)}< L \|\xi - \zeta^i\|.\label{eq:worst-q-r1:L}\end{equation}

By using the total law of probability, $\sup_{\mathbb P \in \mathcal P^{1}_{N,\epsilon}}\mathbb E_{\mathbb P}[Z(x,\tilde \xi)]$ can be represented as 
{\tiny\begin{align}\sup_{\mathbb P \in \mathcal P^{1}_{N,\epsilon}}\mathbb E_{\mathbb P}[Z(x,\tilde \xi)]=\sup_{\mathbb Q_i \in \mathcal P(\Xi),i \in [N]} \ & \frac{1}{N}\sum_{i \in [N]}\int_{\Xi}Z(x,\xi)d\mathbb Q_{i}(\xi) : \frac{1}{N}\sum_{i \in [N]}\int_{\Xi}\|\xi - \zeta^i\|d\mathbb Q_{i}(\xi) \le \epsilon,
\label{eq:worst-case-exp:total}
\end{align}}
where $\mathbb Q_i$ denotes the transportation plan for moving the $\frac{1}{N}$ mass at $\zeta^i$ for each $i \in [N]$. Note that for any feasible $\mathbb Q := \frac{1}{N}\sum_{i \in [N]}\mathbb Q_i \in \mathcal P^{1}_{N,\epsilon}$, the objective function in \eqref{eq:worst-case-exp:total} is evaluated at
\(\frac{1}{N}\sum_{i \in [N]}Z(x,\zeta^i) + \frac{1}{N}\sum_{i \in [N]}\int_{\Xi}Z(x,\xi) - Z(x,\zeta^i)d\mathbb Q_{i}(\xi) < \frac{1}{N}\sum_{i \in [N]}Z(x,\zeta^i) + \frac{1}{N}\sum_{i \in [N]}\int_{\Xi}L \|\xi - \zeta^i\| d\mathbb Q_{i}(\xi) 
\le \frac{1}{N}\sum_{i \in [N]}Z(x,\zeta^i) + \epsilon L,
\)
where the strict inequality follows from \eqref{eq:worst-q-r1:L} and the last inequality from $\mathbb Q \in \mathcal P^{1}_{N,\epsilon}$. Therefore, $\sup_{\mathbb P \in \mathcal P^{1}_{N,\epsilon}}\mathbb E_{\mathbb P}[Z(x,\tilde \xi)]$ cannot be achieved by any feasible $\mathbb Q \in \mathcal P^{1}_{N,\epsilon}$ in this case. Instead, $\sup_{\mathbb P \in \mathcal P^{1}_{N,\epsilon}}\mathbb E_{\mathbb P}[Z(x,\tilde \xi)]$ is asymptotically achieved by $\mathbb Q_{x,\Delta}$, defined in \eqref{eq:dist:r1:R:2}, as $\Delta$ ever approaches to zero. 
It is easy to see that $\mathbb Q_{x,\Delta}$ is in $\mathcal P^{1}_{N,\epsilon}$ for any $\Delta \in (0,\frac{1}{N}]$ since $\|d_{\pi^*}\| = 1$. Note that, for $\epsilon > 0$, the expectation is evaluated in the limit as follows: $\lim_{\Delta \rightarrow 0^+}\frac{1}{N}\sum_{i \in [N] \setminus \{i'\}}Z(x,\zeta^i) + (\frac{1}{N}-\Delta)Z(x,\zeta^{i'}) + \Delta Z(x,\zeta^{i'}+ \frac{\epsilon}{\Delta}d_{\pi^*}) = \frac{1}{N}\sum_{i \in [N]}Z(x,\zeta^i) + \lim_{\Delta \rightarrow 0^+}\Delta (Z(x,\zeta^{i'}+ \frac{\epsilon}{\Delta}d_{\pi^*})-Z(x,\zeta^{i'}))$.
For any $\hat \pi \in \Pi_{i',x}^*$
, the limit term satisfies
\(
\lim_{\Delta \rightarrow 0^+}\Delta (Z(x,\zeta^{i'}+ \frac{\epsilon}{\Delta}d_{\pi^*})-Z(x,\zeta^{i'})) = \lim_{\Delta \rightarrow 0^+} \Delta \left[\max_{\pi \in \Pi} \pi^T (Tx + \zeta^{i'} + \frac{\epsilon}{\Delta}d_{\pi^*}) - \hat\pi^T (Tx + \zeta^{i'})\right]
= \lim_{\Delta \rightarrow 0^+} \left[\max_{\pi \in \Pi}  {\epsilon}\pi^\top d_{\pi^*} + \Delta(\pi - \hat\pi)^\top (Tx + \zeta^{i'})\right] = \epsilon L,
\)
where the last equality is obtained by interchanging the limit and the maximum; this is valid because $f_\Delta(\pi) := {\epsilon}\pi^\top d_{\pi^*}+ \Delta(\pi - \hat\pi)^\top (Tx + \zeta^{i'})$ converges uniformly to $f_0(\pi)$ over compact $\Pi$. 
Therefore, we have $\lim_{\Delta \rightarrow 0^+} \max_{\pi \in \Pi} f_\Delta(\pi) = \max_{\pi \in \Pi}\lim_{\Delta \rightarrow 0^+}  f_\Delta(\pi)= \max_{\pi \in \Pi}f_0(\pi) = \epsilon\max_{\pi\in \Pi} {\pi}^\top d_{\pi^*} = \epsilon L$, since an upper bound
$\max_{\pi\in \Pi} {\pi}^\top d_{\pi^*} \le \max_{\pi\in \Pi} \|\pi\|_q \|d_{\pi^*}\| = \max_{\pi\in \Pi} \|\pi\|_q = L$
is achieved by $\pi = \pi^*$.
\qed

\section{Proof of Proposition \ref{prop:2-wass-worst-dist}}
\label{appendix:proof:2-wass-worst-dist}
First, note that \(\lambda^* > 0\) because if \(\lambda = 0\), the dual problem becomes \(\sup_{\xi \in \Xi} Z(x, \xi) = \sup_{\xi \in \Xi, \pi \in \Pi} \pi^\top Tx + \pi^\top \xi = \infty\), as we can set \(\xi = \alpha \pi\) for \(\alpha \to \infty\) with any \(\pi \in \Pi\). Therefore, without loss of generality, we have
{\tiny\begin{align}
    \sup_{\mathbb{P} \in \mathcal{P}^{2}_{N,\epsilon}}
    \mathbb{E}_{\mathbb{P}}[Z(x, \tilde{\xi})] &= \min_{\lambda > 0} \epsilon^2 \lambda + \frac{1}{N} \sum_{i \in [N]} \sup_{\xi \in \Xi, \pi \in \Pi} \pi^\top (Tx + \xi) - \lambda \|\xi - \zeta^i\|^2 \nonumber \\
    &= \min_{\lambda > 0} \epsilon^2 \lambda + \frac{1}{N} \sum_{i \in [N]} \sup_{\pi \in \Pi} \pi^\top Tx + \sup_{\xi \in \Xi} \pi^\top \xi - \lambda \|\xi - \zeta^i\|^2.\label{eq:worst-case-dist:r2:inner-supremum}
\end{align}}
Note first that, when $\Xi = \mathbb R^k_+$, $\sup_{\xi \in \Xi} \pi^\top \xi - \lambda \|\xi - \zeta^i\|^2 = \sup_{\xi \in \mathbb R^k} \pi^\top \xi - \lambda \|\xi - \zeta^i\|^2$ since $\pi \ge 0$ and $\zeta^i \ge 0$. Therefore, for both cases of $\Xi = \mathbb R^k$ and $\mathbb R^k_+$, the innermost supremum problem in \eqref{eq:worst-case-dist:r2:inner-supremum} corresponds to the conjugate function of \( f(\xi) := \lambda\|\xi - \zeta^i\|^2 = \frac{1}{2}\| \sqrt{2\lambda}I \xi - \sqrt{2\lambda}\zeta^i \|^2 \). The conjugate function can be derived using the conjugate of the composition of a squared norm with an affine transformation (see Example 3.27 and p.95 of \cite{boyd2004convex}):
\(
f^*(\pi) = \frac{1}{2}\left\|\frac{1}{\sqrt{2\lambda}}\pi\right\|_q^2 + \sqrt{2\lambda}(\zeta^i)^\top \frac{1}{\sqrt{2\lambda}} \pi =  \pi^\top \zeta^i + \frac{1}{4\lambda}\|\pi\|_q^2.
\)
This represents \( \sup_{\mathbb P \in \mathcal P^{2}_{N,\epsilon}} \mathbb E_{\mathbb P}[Z(x,\tilde \xi)] \) as follows: 
\begin{equation}\tiny
\sup_{\mathbb{P} \in \mathcal{P}^{2}_{N,\epsilon}} \mathbb{E}_{\mathbb{P}}[Z(x, \tilde{\xi})] = \min_{\lambda > 0} \epsilon^2 \lambda + \frac{1}{N} \sum_{i \in [N]} \sup_{\pi \in \Pi} \pi^\top (Tx + \zeta^i) + \frac{1}{4\lambda} \|\pi\|_q^2.
\label{eq:worst-case-dist:r2:pi}
\end{equation}
Now, define \( f_i(\lambda) := \sup_{\pi \in \Pi} \pi^\top (Tx + \zeta^i) + \frac{1}{4\lambda} \|\pi\|_q^2 \), which is a convex and nonincreasing function of \(\lambda \in \mathbb{R}_{++}\). Therefore, \(\lambda^*>0\) minimizes \(\epsilon^2 \lambda + \frac{1}{N} \sum_{i \in [N]} f_i(\lambda)\) if and only if 
\(
0 \in \epsilon^2 + \frac{1}{N} \sum_{i \in [N]} \partial f_i(\lambda^*), 
\)
where \(\partial f_i(\lambda^*) = \text{conv}\left\{ -\frac{1}{4(\lambda^*)^2} \|\pi\|_q^2 : \pi \in \Pi^*_{i, x, \lambda^*} \right\}\), meaning
\(\frac{1}{N} \sum_{i \in [N]} \sum_{\pi \in \Pi^*_{i, x, \lambda^*}} \mu^i_{\pi} \frac{\|\pi\|_q^2}{4(\lambda^*)^2} = \epsilon^2, \)
for some \((\mu^i_{\pi})_{\pi \in \Pi^*_{i, x, \lambda^*}} \geq 0\) summing to one for each \( i \in [N] \). It follows from \eqref{eq:worst-case-prob-dist:r2:dual} that \(\mathbb{Q}_x\), as defined in \eqref{eq:worst-prob-dist:r2} with the \((\mu^i_{\pi})\)’s satisfying \eqref{eq:worst-case-prob-dist:r2:dual}, lies in \(\mathcal{P}^{2}_{N,\epsilon}\) and is optimal since
\begin{subequations}
    \tiny\begin{align}
        \mathbb{E}_{\mathbb{Q}_x}[Z(x, \tilde{\xi})] &= \epsilon^2 \lambda^* + \frac{1}{N} \sum_{i \in [N]} \sum_{\pi \in \Pi^*_{i, x, \lambda^*}} \mu^i_{\pi} Z(x, \zeta^i + \frac{\|\pi\|_q}{2\lambda^*}d_\pi) - \epsilon^2 \lambda^* \\
        &= \epsilon^2 \lambda^* + \frac{1}{N} \sum_{i \in [N]} \sum_{\pi \in \Pi^*_{i, x, \lambda^*}} \mu^i_{\pi} \left( \max_{\pi' \in \Pi} (\pi')^\top (Tx + \zeta^i + \frac{\|\pi\|_q}{2\lambda^*}d_\pi) - \frac{\|\pi\|_q^2}{4\lambda^*} \right) \label{eq:worst-prob-dist:r2:1} \\
        &\geq \epsilon^2 \lambda^* + \frac{1}{N} \sum_{i \in [N]} \sum_{\pi \in \Pi^*_{i, x, \lambda^*}} \mu^i_{\pi} \left(\pi^\top (Tx + \zeta^i) + \frac{\|\pi\|^2_q}{4\lambda^*} \right)\label{eq:worst-prob-dist:r2:2} \\
        &= \epsilon^2 \lambda^* + \frac{1}{N}\sum_{i \in [N]}f_i(\lambda^*) = \sup_{\mathbb{P} \in \mathcal{P}^{2}_{N,\epsilon}} \mathbb{E}_{\mathbb{P}}[Z(x, \tilde{\xi})]. \label{eq:worst-prob-dist:r2:3}
    \end{align}
\end{subequations}
Here, \eqref{eq:worst-prob-dist:r2:1} follows from \eqref{eq:worst-case-prob-dist:r2:dual}, \eqref{eq:worst-prob-dist:r2:2} from the suboptimality of \(\pi\) for the inner maximization and $\pi^\top d_\pi = \|\pi\|_q$, and \eqref{eq:worst-prob-dist:r2:3} from the definition of \( f_i(\lambda^*)\) and that $\sum_{\pi \in \Pi^*_{i,x,\lambda^*}}\mu^i_\pi = 1, \forall i \in [N]$.
\qed

\section{Proof of Proposition \ref{prop:newsvendor}}\label{appendix:proof:newsvendor}
Proposition \ref{prop:r1:worst-expectation} indicates that $x^{1}_{N}(\epsilon) = x^{\texttt{SAA}}_1 = -\zeta = \arg\min_{x\ge 0} (\texttt c - \texttt p) x + Z(x,\zeta)$.
When $r=2$, according to \eqref{eq:worst-case-dist:r2:inner-supremum}, the worst-case expectation problem can be expressed as
\begin{equation}\tiny
\min_{\lambda > 0} \, \epsilon^2 \lambda + \max_{\pi \in [0, \texttt{p}]} \left( \pi (x + \zeta) + \frac{1}{4\lambda} \pi^2 \right) = \min_{\lambda > 0} \, \epsilon^2 \lambda + \max_{\pi \in \{0, \texttt{p}\}} \left( \pi (x + \zeta) + \frac{1}{4\lambda} \pi^2 \right),
\label{eq:worst-case:newsvendor}
\end{equation}
where the equality holds because the inner problem is a convex maximization, meaning that an optimum occurs at a boundary point. 
Since \eqref{eq:worst-case:newsvendor} is a convex minimization over \(\lambda > 0\), \(\lambda^*\) is optimal if and only if \(0 \in \partial f(\lambda)\), where  
\begin{equation*}\tiny
\partial f(\lambda) =
\begin{cases}
\epsilon^2, \ & \mbox{ for } \lambda: \texttt{p}(x+\zeta) + \frac{1}{4\lambda}\texttt{p}^2 < 0,\\
\epsilon^2 - \frac{1}{4\lambda^2}\texttt{p}^2, \ &  \mbox{ for } \lambda: \texttt{p}(x+\zeta) + \frac{1}{4\lambda}\texttt{p}^2 > 0,\\
\left\{\epsilon^2 - \frac{\mu}{4\lambda^2}\texttt{p}^2: \mu \in [0,1]\right\}, \ & \mbox{ for } \lambda: \texttt{p}(x+\zeta) + \frac{1}{4\lambda}\texttt{p}^2 = 0,
\end{cases}\end{equation*}
Note that for \(x > -\zeta - \frac{\epsilon}{2}\), 
\(\lambda^* = \frac{\texttt{p}}{2\epsilon}\) and \(\Pi^*_{x,\lambda^*} = \{\texttt{p}\}\) since 
$\texttt{p}(x+\zeta) + \frac{1}{4\lambda^*}\texttt{p}^2 > 0$ and $\partial f(\lambda^*) = \epsilon^2 - \frac{\texttt{p}^2}{4\lambda^2} = 0$.
Therefore, 
\(
\sup_{\mathbb{P} \in \mathcal{P}^{2}_{N,\epsilon}} \mathbb{E}_{\mathbb{P}}[Z(x, \tilde{\xi})] = \texttt{p} (x + \zeta) + \epsilon \texttt{p} \mbox{ for } x \ge -\zeta - \frac{\epsilon}{2}.
\)
Proposition \ref{prop:2-wass-worst-dist} indicates that the worst-case distribution is \( \mathbb Q_x=\delta_{\zeta + \frac{\|\texttt{p}\|_q}{2\lambda^*} d_{\texttt{p}}} = \delta_{\zeta + \epsilon} \), corresponding to a scenario where the demand is realized as \( -(\zeta + \epsilon) \). 



Otherwise, i.e., for \(x \leq -\zeta - \frac{\epsilon}{2}\), we have \(\lambda^* = \frac{\texttt{p}}{-4(x+\zeta)}\) and \(\Pi^*_{x,\lambda^*} = \{0, \texttt{p}\}\), as it satisfies $0 \in \partial f(\lambda^*)$ 
with \(\mu = \frac{\epsilon^2}{4(x+\zeta)^2} \leq 1\); specifically, $\epsilon^2 - \frac{\mu \texttt{p}^2}{4 (\lambda^*)^2}  = \epsilon^2 - 4\mu(x+\zeta)^2 = 0$.
Proposition \ref{prop:2-wass-worst-dist} indicates that the worst-case probability distribution is supported at \(\zeta + \frac{\texttt{p}}{2\lambda^*} = -2x-\zeta\) and \(\zeta\) and is given by  
\(
\mathbb Q_x = \mu \delta_{-2x-\zeta} + \left(1-\mu\right) \delta_{\zeta},
\)  
where \(\mu = \frac{\epsilon^2}{4(x+\zeta)^2}\). Therefore,  
\(
\sup_{\mathbb{P} \in \mathcal{P}^{2}_{N,\epsilon}} \mathbb{E}_{\mathbb{P}}[Z(x, \tilde{\xi})] = \texttt{p} \frac{\epsilon^2}{4(x+\zeta)^2} \left(x - (2x+\zeta)\right) = -\frac{\epsilon^2 \texttt{p}}{4(x+\zeta)} \ \text{for } x \leq -\zeta - \frac{\epsilon}{2}.
\)

\begin{figure}[t!]
\centering
\includegraphics[width=0.33\textwidth]{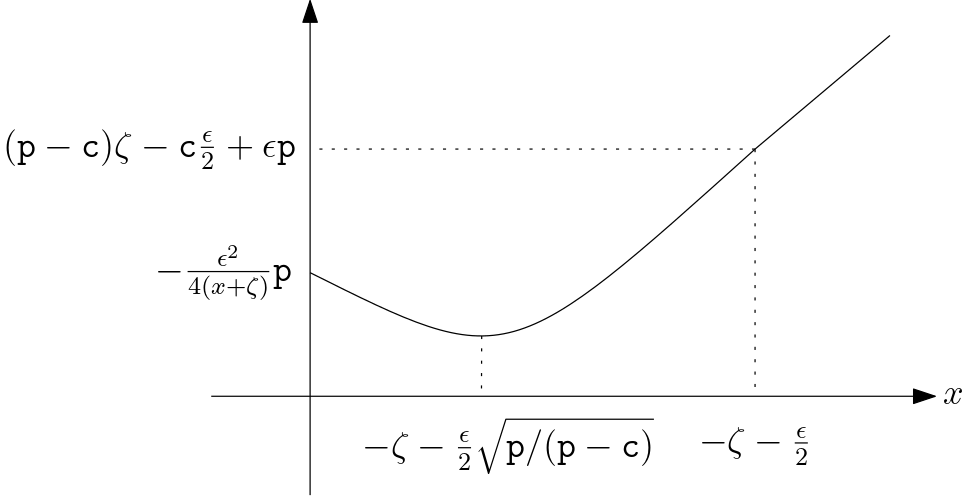}
        \caption{The graph of $g$ when $-\zeta - \frac{\epsilon}{2}\sqrt{\frac{\texttt p}{\texttt p - \texttt c}}>0$.}
        \label{fig:g}
\end{figure}

Therefore, \eqref{prob:newsvendor} becomes $\min_{x \ge 0}g(x)$, where
\begin{align*}\tiny
g(x) := \begin{cases}(\texttt{c} - \texttt{p}) x + \texttt{p}(x+\zeta)+ \epsilon\texttt p=\texttt{c}x + \texttt{p}\zeta+ \epsilon\texttt p,  &\mbox{ for } x > -\zeta-\frac{\epsilon}{2} \\
(\texttt{c} - \texttt{p}) x - \frac{\epsilon^2\texttt{p}}{4(x+\zeta)} 
, &\mbox{ for } x \le -\zeta-\frac{\epsilon}{2}\\
\end{cases}
\end{align*}
which is visualized in Figure \ref{fig:g}. Therefore, the realized solution for a given $\zeta$ is
\(
x^{2}_1(\epsilon) = \max\left\{0, -\zeta - \frac{\epsilon}{2}\sqrt{\frac{\texttt{p}}{\texttt{p} - \texttt{c}}}\right\},
\)  
which depends on the risk appetite \(\epsilon\) and the critical fractile \(\frac{\texttt{p} - \texttt{c}}{\texttt{p}}\).

\bibliographystyle{plain}
\bibliography{reference}

\end{document}